\documentclass[conference]{IEEEtran}
\usepackage{amsmath,amssymb}

\begin{document}
\begin{center}
\large\textbf{Distance Graphs of  Metric Spaces with Rosenbloom -
Tsfasman metric }
\end{center}
\begin{center}{W. B. Vasantha and R. Rajkumar}\\
{Department of Mathematics }\\
{Indian Institute of Technology
Madras, Chennai-600 036.}\\
{email: vasantha , rajkumarr @iitm.ac.in}\\

\end{center}

\vspace{.5 cm}

 \begin{abstract} Rosenbloom and Tsfasman introduced a new
 metric (RT metric) which is a generalization of the Hamming metric.
  In this paper we study the distance graphs of
 spaces $Z_q^n$ and $S_n$ with Rosenbloom -Tsfasman metric. We also describe the degrees of vertices,
 components and  the chromatic number of these graphs.
  Distance graphs of general direct product spaces also described.    \\
\end{abstract}
%

\section{Introduction}

 The study of metric spaces by
using its corresponding distance graphs tells more about its
structure. A lot of research has been done for finding the structure
of the distance graphs on spaces $R^n$ and $Z^n$ with several
metrics like Euclidean, $l_1$ and $l_p$. A particular interest in
this area is to find the chromatic number of these distance graphs.
[ 2,3,4,5,6] give some references in this direction.


 In this paper we mainly focus on the distance graphs of two families of spaces
 $Z_q^n$ and $S_n$ with RT metric given by Rosenbloom and Tsfasman [7] in
 1997. In particular, we investigate the degrees of vertices, components
 and the chromatic numbers for these graphs.\\

\section{ Notations and  Preliminaries}
%
%
%

For any connected graph $G$ we write $nG$ for the graph with $n$
components each isomorphic to $G$. In this paper we denote by
$[G]^m$, the join of the graph $G$ with itself $m$ times.
\begin{center}
$i.e.,~~~~[G]^m = G+ \cdots + G $ ($m$ times).
\end{center}

A complete graph with $n$ vertices is denoted by $K_n$ and a
complete $r$-partite graph $K_{m, \ldots, m}$ is denoted by
$K_r(m)$. The \textit{chromatic number}  of a graph $G$ is denoted
by $\chi (G)$


The distance graph of a metric space is defined in [2,3] as follows.
Suppose $S$ is a subset of a metric space $X$ with metric $d$.
Denote by dist$(S)$ the set of all distances realized by pairs of
distinct points in $S$, i.e.,

\begin{center}
 dist$(S) = \{d(x, y) : x,y \in S$ and $x\neq y\}$.
\end{center}

For each subset $D$ of dist$(S)$ the \textit{distance graph}
$G(S,D)$ with \textit{distance set} $D$ is the graph with vertex set
$S$ and edge set $E(S,D)=\{ xy:x,y \in S$ and $d(x,y) \in D \}$.\\

Thus for a metric space $X$, $\displaystyle \bigcup_{{k \in}dist(X)}
G(X,\{k\})$ is a complete graph on $X$.\\


 Here we give the definition of the RT metric on  $Z_q^n$.[1]

 Let $Z_q=\{ 0,1,\ldots ,q-1\}$. Let $x = (x_1, \ldots, x_n)\in
 Z_q^n$. Then
\begin{center}
 $$ \omega(x)=\begin{cases} max~ \{i ~| ~x_i \neq 0 \}~,~x \neq0
\\ 0  ~~~~~~~~~~~~~~~~~~~~~,~x=0 \end{cases}$$
\end{center}

 is called the \textit{RT weight} of $x$.

 The \textit{RT distance} between $x$ and $y$ in
$Z_q^n$ is defined by

\begin{center}
$\rho(x,y)=\omega(x-y).$
\end{center}

   $\rho$ is a metric on $Z_q^n$.

%


\vspace{.5 cm}

\section { The distance graphs on $Z_q^n$ with RT metric}

In this section we obtain  the distance graph $G(Z_q^n,D)$ of
$(Z_q^n,\rho)$. At the end of this section we generalize this to
general direct product space.  Note that for $(Z_q^n,\rho)$,
dist$(Z_q^n) =
\{1,\ldots ,n \}$.\\


%
%

%
%


\noindent \textit{\textbf{Theorem 3.1:} For the metric space
$(Z_q^n,\rho)$, where $\rho$ denotes the RT metric, suppose
$D=\{d_1,\ldots ,d_k\}\subseteq \{1, \ldots , n\}$ with $ 1 \leq d_1
< \cdots < d_k \leq n $ then the distance graph $G (Z_q^n, D)$ of
$(Z_q^n,\rho)$ is isomorphic to}

\textit{\begin{center} $\displaystyle
q^{n-d_k}\left[q^{d_k-d_{k-1}-1}\left[ \cdots \left[q^{d_2- d_1 -1}
K_q(q^{d_1-1})\right]^q \cdots \right]^q \right]^q .$
\end{center}}
\vspace{.5 cm}

 \noindent \textbf{\textit{Proof:}} Given the distance set
 $D=\{d_1,\ldots ,d_k\}\subseteq \{1, \ldots ,n\}$ with $ 1 \leq d_1
< \cdots < d_k \leq n $.

Clearly $G(Z_q^n,D)= \displaystyle \bigcup
^k_{i=1}G(Z_q^n,\{d_i\})$.\\

 First we give the structure of the graphs $G(Z_q^n,\{d\})$ for $d
 \in \{1, \ldots , n\}$.

Let $x = (x_1, \ldots , x_d,\ldots , x_n)\in Z_q^n$ be a fixed
vector.

 A vector $y = (y_1, \ldots , y_d,\ldots , y_n) \in
 Z_q^n$ is adjacent with $x$  in $G (Z_q^n, \{d\})$  if and only if $\rho (x,y)=d$. But
 $\rho (x,y)=d \Leftrightarrow x_d \neq y_d $ and $
 x_i = y_i ,d<i \leq n$.

 By fixing the coordinate $y_d$ in $y$ such that $x_d \neq y_d$, the
 remaining coordinates $y_1, \ldots , y_{d-1}$ can be chosen in
 $q^{d-1}$ different ways. Thus we get a set of  $q^{d-1}$ vectors
 such that all of its $d^{th}$ coordinates are same and equal to $y_d$ and they are at
 a  distance $d$ form $x$. We also note that this set of
 $q^{d-1}$ vectors are not adjacent with each other since the
 distance between any two of them is strictly less than $d$.

 Now the $d^{th}$ coordinate in $y$ can vary in $(q-1)$ different
 ways. By including the vector $x$ , we get totally $q$ sets of
 $q^{d-1}$ vectors such that the vectors are not adjacent with each other
 with in each $q$ sets, but all $q^{d-1}$ vectors in each one of
 these $q$ sets are adjacent to the remaining sets of vectors. Thus
 they form a complete $q$-partite graph $ K_q(q^{d-1})$.
 Now the vector $x$ can vary in $q^{n-d}$ ways. Thus correspondingly
 we get $q^{n-d}$ copies of $ K_q(q^{d-1})$.
\begin{center}
 i.e., $G (Z_q^n,\{d\} )\cong q^{n-d}K_q(q^{d-1})$.
\end{center}

 So for each elements of distance set $D$ we have
  $G (Z_q^n,\{d_i\} )\cong q^{n-d_i}K_q(q^{d_i-1})$.

Now we construct the graph $G(Z_q^n,\{d_1,d_2\})$.

From the definition of the distance graph of a metric space two
vertices in $G(Z_q^n,\{d_1,d_2\})$ are adjacent if and only if the
RT distance between them is either $d_1$ or $d_2$.

So  $G(Z_q^n,\{d_1,d_2\})$ is the edge disjoint union of
$G(Z_q^n,\{d_1\})$ and $G(Z_q^n,\{d_2\})$. Thus to get the graph
$G(Z_q^n,\{d_1,d_2\})$  it is enough to add edges to
$G(Z_q^n,\{d_2\})$ in such a way that if $(u,v)$ is an edge added to
$G(Z_q^n,\{d_2\})$ then $(u,v)$ is an edge in $G(Z_q^n,\{d_1\})$.

Now the graph $G(Z_q^n,\{d_2\})$ can be written as ,

$G(Z_q^n,\{d_2\})=G_1 \cup \cdots \cup G_{q^{n-d_2}}$, with each
$G_i \cong K_q(q^{d_2-1}),i=1, \ldots , q^{n-d_2}. $

For each $i=1, \ldots , q^{n-d_2} $ let $V_i$ be the vertex set of
$G_i$. Thus the space $Z_q^n$ has been partitioned into $q^{n-d_2}$
disjoint sets $V_i ~ ,i=1, \ldots , q^{n-d_2}$ with $\mid V_i \mid =
q^{d_2}$.

Since each $G_i \cong K_q(q^{d_2-1}) $, its corresponding vertex set
$V_i$ is also partitioned in to disjoint sets $V_{ji}$ such that,

$V_i =  V_{1i} \cup \cdots \cup V_{qi}$ , $\mid V_{ji} \mid
=q^{d_2-1}, ~j = 1, \ldots , q$

 with the condition that each edge in $G_i$ joins a vertex in $V_{li}$
 and a vertex in $V_{ji},~ l,j= 1, \ldots , q $ with $ l \neq j$, for all
 $i=1, \ldots , q^{n-d_2}$.

In $G(Z_q^n,\{d_2\})$ two vertices $x $ and $y$ are adjacent if and
only if $\rho (x,y) = d_2$.

In the construction of $G(Z_q^n,\{d_2\})$  we note that the
partitioned sets of $Z_q^n$ have the following properties:

For any $x = (x_1,\ldots , x_n)$ and $y = (y_1 , \ldots , y_n) \in
Z_q^n$,
\begin{itemize}

\item[(1)] If $x \in V_i$ and $y \in V_j$ for some $i,j=1, \ldots , q^{n-d_2}$ with $i
\neq  j $  then $(x_{d_2+1}, \ldots ,x_n) \neq (y_{d_2+1}, \ldots
,y_n) $; i.e., $\rho (x,y) >d_2$.

\item[(2)] If $x,y \in V_i$ for some $i=1, \ldots , q^{n-d_2}$ then
$x_r = y_r$ $ \forall~ r = d_2+1, \ldots , n$; i.e., $\rho (x,y)
\leq d_2$.

\item[(3)] If $x,y \in V_{ji}$ for some $i=1, \ldots , q^{n-d_2}$
and $j=1, \dots , q$ then $x_r = y_r$ $ \forall~ r = d_2, \ldots ,
n$; i.e., $\rho (x,y) < d_2$.

\item[(4)] If $x \in V_{li}$ and $y \in V_{ji}$ for some $i=1, \dots
, q^{n-d}$ and $l,j = 1, \ldots , q$ with $l \neq j $ then $ x_{d_2}
\neq y_{d_2}$ and $x_r = y_r$ $ \forall~ r = d_2+1, \ldots , n ; $
i.e., $\rho (x,y)= d_2$.

\end{itemize}

For a fixed $i \in \{1, \ldots , q^{n-d_2}\}$ consider the component
$G_i$ of $G(Z_q^n,\{d_2\})$.

Let $x= (x_1,\ldots , x_n) \in V_i$ then $x \in V_{ji}$ for some
$j\in \{1, \ldots ,q\}$ (say). Consider $y = (y_1 , \ldots , y_n)
\in Z_q^n$ with $\rho(x,y)=d_1$. We show that any such $y$ must be
in $ V_{ji}$. Suppose

\begin{itemize}
\item[(i)]  if $y \in V_k$ for some $k \neq i$, $k=1, \ldots , q^{n-d}$
then by using property $(1)$ and $d_1 <d_2$ , we have $\rho (x ,y) >
d_2 > d_1$ , which is a contradiction. So $y \in V_i$.

\item[(ii)] if $y \in V_{li}$ for some $l\neq j, ~l=1, \ldots , q $
then by using property $(4)$ we have $ \rho (x, y)= d_2$, which is a
contradiction. So $y \in V_{ji}$.

\end{itemize}

Hence all the vectors which are at a distance $d_1$ from  $x$ must
be in $V_{ji}$ .

Since $x$ is an arbitrary point in $V_{ji}$ and the component of
$G(Z_q^n,\{d_1\})$ which has $x$ as a vertex is isomorphic to
$K_q(q^{d_1-1})$, it follows that all the edges in the component
containing $x$ must be added in $V_{ji}$.

Since each component $K_q(q^{d_1-1})$ of $G(Z_q^n,\{d_1\})$ has
$q^{d_1}$ vertices and $| V_{ji} | =q^{d_2-1}$, we can attach
$q^{d_2-d_1-1}$ number of components $K_q(q^{d_1-1})$ in $V_{ji}$.

In this way attach components $K_q(q^{d_1-1})$ of $G(Z_q^n,\{d_1\})$
in all partitioned vertex sets of $G_i$.

Thus after adding edges from $G(Z_q^n,\{d_1\})$ to the component
$G_i$ of $G(Z_q^n,\{d_1,d_2\})$ which we are considering is
isomorphic to $\left[ q^{d_2-d_1-1} K_q(q^{d_1-1})\right]^q$.

One can in this way attach $q^{d_2-d_1}$ components $K_q(q^{d_1-1})$
of $G(Z_q^n,\{d_1\})$ in each component of $G(Z_q^n,\{d_2\})$.

  But the number of components of $G(Z_q^n,\{d_2\})$ is $q^{n-d_2}$,
  thus totally we are attaching $q^{d_2-d_1}q^{n-d_2}=q^{n-d_1}$
  i.e, all the components of $G(Z_q^n,\{d_1\})$ in $G(Z_q^n,\{d_2\})$.

and hence
\begin{center}
 $G\left(Z_q^n,\{d_1,d_2\}\right) \cong q^{n-d_2}
 \left[q^{d_2-d_1-1}K_q(q^{d_1-1})\right]^q$.
\end{center}

Now by using the similar type of argument we can construct
$G(Z_q^n,\{d_1,d_2,d_3\})$.

 Note that
$G(Z_q^n,\{d_1,d_2,d_3\}) = G(Z_q^n,\{d_1,d_2\}) \cup
G(Z_q^n,\{d_3\})$.

Any component of $G(Z_q^n,\{d_1,d_2\})$ is isomorphic to $\left[
q^{d_2-d_1-1} K_q(q^{d_1-1})\right]^q$. By attaching these
components in $G(Z_q^n,\{d_3\})$  as in the method described in the
previous case one can show that any component of
$G(Z_q^n,\{d_1,d_2,d_3\})$ is isomorphic to
\begin{center}
$\displaystyle \left[q^{d_3-d_2-1}\left[ q^{d_2-d_1-1}
K_q(q^{d_1-1})\right]^q \right]^q$.
\end{center}
Since $G(Z_q^n,\{d_3\})$ has $q^{n-d_3}$ components, we have
\begin{center}
$G(Z_q^n,\{d_1,d_2,d_3\}) \cong q^{n-d_3}\displaystyle
\left[q^{d_3-d_2-1}\left[ q^{d_2-d_1-1} K_q(q^{d_1-1})\right]^q
\right]^q$.

\end{center}

Proceeding in this way, at the $k^{th}$ stage, we can arrive at the
graph $G(Z_q^n,D)$ as described in the statement of this theorem.
\hfill{$\Box$}\\

%
%
%
%
%


\noindent\textit{\textbf{Note 3.1: }}Let $p, r_i ,s_i$
$(i=1,\ldots,k)$ be positive integers with $r_1 < \cdots <r_k$ and
$s_1 < \cdots < s_k$ and $p>1$ then by using induction it is easy to
show that $p^{r_1}+ \cdots + p^{r_k}= p^{s_1}+ \cdots + p^{s_k}$ if
and only if $r_i = s_i$, $\forall i=1,\ldots,k$.\\


\noindent\textit{\textbf{Corollary 3.1: } For the metric space
$(Z_q^n,\rho)$, where $\rho$ denotes the RT metric, suppose
$D=\{d_1,\ldots ,d_k\}\subseteq \{1, \ldots , n\}$  with $ 1 \leq
d_1 < \cdots < d_k \leq n $ then the following holds.
\begin{itemize}
\item[(1)] $G(Z_q^n , D)$ is regular of degree $(q-1)\displaystyle
\sum _{i=1}^{k} q^{d_i-1}$.
\item[(2)] $G(Z_q^n, D)$ is connected if and only if $n \in D$.
\item[(3)] The components of $G(Z_q^n,D)$ and $G(Z_q^m , D_1)$
are isomorphic if and only if $D = D_1$.
\item[(4)] $\chi \left(G (Z_q^n,D )\right)=q^k$.
\item[(5)] For any two distance sets $D_1$ and $D_2$,
 $\chi\left( G (Z_q^n,D_1 )\right) = \chi \left( G (Z_q^n,D_2 )\right)$
 if and only if  $|D_1|=|D_2|$.
\end{itemize}}

 \noindent\textit{\textbf{Proof:}}
\noindent(1) Let $x \in Z_q^n$. Then $\displaystyle \mid \{y \in
Z_q^n : \rho (x,y)= d\}\mid  = (q-1)q^{d-1}.$

Hence the degree of a vertex in $G(Z_q^n , D)$ is $\displaystyle
\mid \{y \in Z_q^n : \rho (x,y)\in D \}\mid = (q-1)\displaystyle
\sum _{i=1}^{n} q^{d_i-1}.$\\

\noindent(2) By theorem (3.1), $G(Z_q^n, D)$ is connected if and
only if $q^{n-d_k}=1$ ; if and only if  $d_k =n $.\\

\noindent(3) If $D = D_1$ then by theorem (3.1),
 $G(Z_q^n, D)$ and $G(Z_q^m, D)$  differ only by the number of
 components.

 Suppose $D_1=\{d_1^{'},\ldots , d_r^{'}\} \subseteq \{1,\ldots , m\}$ and
 a  component $H$  of  $G(Z_q^n, D)$ is  isomorphic to a component
 $H^{'}$ of  $G(Z_q^m, D_1)$ then we have to show that $D = D_1$.
 Since $H \cong H^{'}$  the number of vertices and the degrees of vertices
 in these components
 are same. Thus from theorem (3.1), we have $\mid D \mid = \mid D_1 \mid$ and so
\begin{center}
 $(q-1)\displaystyle
\sum _{i=1}^{k} q^{d_i-1}=(q-1)\displaystyle \sum _{i=1}^{k}
q^{d_i^{'}-1}$ .
\end{center}

 Now by using note $3.1$
we must have $d_i =d_i^{'}$, $\forall~ i=1, \ldots ,k$ and hence $D = D_1$.  \\

\noindent(4) Using theorem (2.1),(3.1) and  $\chi(K_q(q^{d-1}))=q$,
the result follows.\\

\noindent(5) Follows directly from (4). \hfill{$\Box$}\\

\subsection{Distance graphs of general direct product spaces}
 Consider the collection of non-empty sets
$\{X_i\}_{i=1}^n$ with $|X_i|=q_i$, $\forall ~ i=1,\ldots ,n$. Let
$X := \displaystyle \prod _{i=1}^n X_i$

We define the RT metric on $X$ as follows:

For $x=(x_1, \ldots , x_n)$ and $y=(y_1, \ldots , y_n)$ in $X$,
\begin{center}
 $$ \rho (x,y)=\begin{cases} max~ \{i ~| ~x_i \neq y_i \},  1\leq i \leq n
\\ 0  ~~~~~~~~~~~~~~~~~~~~~,~otherwise \end{cases}$$
\end{center}
\vspace{.5 cm}
 Then we can get the distance graph $G(X,D)$ of the
matric space $(X, \rho)$ , by the same method as described in
theorem 3.1. We state this as follows.

\noindent \textit{\textbf{Theorem 3.2:} Consider $X := \displaystyle
\prod _{i=1}^n X_i$ , where  $\{X_i\}_{i=1}^n$ is a collection of
nonempty sets with $|X_i|=q_i$ ,$\forall ~ i=1,\ldots ,n$. Then for
the metric space $(X,\rho)$, where $\rho$ denotes the RT metric,
suppose $D=\{d_1,\ldots ,d_k\}\subseteq \{1, \ldots , n\}$ with $ 1
\leq d_1 < \cdots < d_k \leq n $ then the distance graph $G (X, D)$
of $(X,\rho)$ is isomorphic to}

{\footnotesize \begin{eqnarray} \nonumber \displaystyle && q_{d_k+1}
\cdots q_n \left[q_{d_{k-1}+1} \cdots q_{d_k -1}\left[ \cdots
\right.\right.\\\nonumber
 &&\cdots \left. \left.
\left[q_{d_1+1} \cdots q_{d_2-1} K_{q_{d_1}}(q_1q_2 \cdots
q_{d_1-1})\right]^{q_{d_2}} \cdots \right]^{q_{d_{k-1}}}
\right]^{q_{d_k}}.\nonumber
\end{eqnarray}}

\vspace{.5 cm}

 \noindent \textsc{Remark 3.1:} In theorem 3.2 if
$d_{i+1}-d_i=1$ for some $i=1,\ldots, n-1$ then  assume that the
factor $q_{d_{i-1}+1} \cdots q_{d_i -1}$ is equal to 1 and if $d_1
=1$ then also assume
that the factor  $q_1q_2 \cdots q_{d_1-1}$ equals 1.\\

 \noindent \textsc{Remark 3.2:} The results
in Corollory 3.1 can be easily modified to this general direct
product spaces.\\





%
%
%
%
%
%
%
%
%

 \section{The distance graphs on $S_n$ with RT
 metric}


 We define the RT distance on symmetric group $S_n$ on the set
 $W=\{1, \ldots ,n\}$ as follows:

 Take $\alpha \in S_n$. Let $W(\alpha)=\{i \in W : \alpha (i) \neq i\}$
 and $\omega(\alpha)=$ max $\{i \in W: \alpha (i) \neq i\}$.
$\omega(\alpha)$ is called the RT weight of $\alpha$.

 Now we define
the RT distance between $\alpha$ and $\beta$ in $S_n$ by
 $\rho (\alpha , \beta)=$ $
 \omega (\alpha^{-1} \beta ) $.\\

%

 \noindent\textit{\textbf{Theorem 4.1:}
 $(S_n,\rho)$ is a metric space,
 where $\rho$ denotes the RT distance.}


\noindent\textit{\textbf{Proof :}} We prove that $\rho$ satisfies
all the axioms of a metric.

 $(i)$ Clearly $\rho (\alpha , \beta)=\rho (\beta, \alpha )$, for
 all $\alpha , \beta \in S_n$.

 $(ii)$ Since the identity permutation is the only one that
 keeps every element of $W$ fixed, we have   $\rho (\alpha , \beta)\geq 0$
and    $\rho (\alpha , \beta)=0$ only when $\alpha =\beta$.

$(iii)$ Next an element of $W$ is moved by $\alpha \beta$ must be
moved by at least one of the permutations $\alpha , \beta$. Hence
$W(\alpha \beta)\subset W(\alpha)\cup W(\beta)$ and consequently,
$\omega (\alpha \beta)\leq \omega (\alpha)+ \omega(\beta)$.

Now let $\lambda,\mu ,\nu$ be any three elements of $S_n$. We have

$\rho (\lambda , \nu )=\omega(\lambda^{-1} \nu   )=\omega(\lambda
 ^{-1}\mu \mu ^{-1}\nu )\leq  \omega(\lambda ^{-1}\mu) + \omega(\mu ^{-1} \nu  )
 =\rho (\lambda , \mu )+\rho (\mu , \nu )$.

 So the function $\rho$ satisfies the triangle inequality. In this
 way we see that $\rho$ is a metric on $S_n$.\hfill{$\Box$}\\

\noindent \textit{\textbf{Remark:}} Since two permutations cannot
differ exactly in the first position, we have dist$(S_n)=\{2,3,4, \ldots,n\}$ for $(S_n,\rho)$.\\

The next theorem gives the structure of the distance graph for
$(S_n,\rho)$.\\




 \noindent \textit{\textbf{Theorem 4.2:} For the metric space
$(S_n,\rho)$, where $\rho$ denotes the RT metric, suppose
$D=\{d_1,\ldots ,d_k\}\subseteq \{2,3,4, \ldots,n\}$ with  $ 2 \leq
d_1 < \cdots < d_k \leq n $ then distance graph $G(S_n,D)$ of
$(S_n,\rho)$ is isomorphic to }
 \textit{\begin{center} $\displaystyle
\frac{n!}{d_k!}\left[\frac{(d_k-1)!}{d_{k-1}!}\left[ \ldots
\left[\frac{(d_2-1)!}{d_1!} K_{d_1}((d_1-1)!)\right]^{d_2} \ldots
 \right]^{d_{k-1}}\right]^{d_k}.$
\end{center}}
%
%
%
%


 \noindent\textit{\textbf{Proof:}} Given the distance set
 $D=\{d_1,\cdots ,d_k\}\subseteq \{2, \ldots ,n\}$ with $ 2 \leq
d_1 < \ldots < d_k \leq n $.

Clearly $G(S_n,D)= \displaystyle \bigcup ^k_{i=1}G(S_n,\{d_i\})$.

 First we give the structure of the graphs $G(S_n,\{d\})$ for $d
 \in \{1, \ldots , n\}$.

Fix a permutation  $\alpha \in S_n$. A permutation $\beta \in
  S_n$ is adjacent with $\alpha$ $G(S_n,\{d\})$ in if and only if
  $\rho (\alpha,\beta)=d$. But
  $\rho (\alpha,\beta)=d \Leftrightarrow \alpha (d) \neq \beta (d) $ and $
  \alpha (i) = \beta (i),~~d<i \leq n$.

  Since $\beta$ is a permutation, each $\beta (i) \in \{1,
  \ldots , n\}$ is distinct. The last $n-d$ coordinates of $\beta$
  are such that $\alpha (i) = \beta (i) ,d<i \leq n$. So by fixing
  $\beta (d)$ such that $\alpha(d)\neq \beta (d)$, the remaining
  $\beta (1), \ldots, \beta (d-1)$ can be chosen in $(d-1)!$
  different ways.

  Thus we get a set of  $(d-1)!$ different permutations $\beta$
  such that all of its $d^{th}$ coordinates are same and they are at
  a  distance $d$ form $\alpha$. Also this set of
  $(d-1)!$ permutations are not adjacent with each other since the
  distance between any two of them is strictly less than $d$.

  Now the $d^{th}$ coordinate in $\beta $ can vary in $(d-1)$ different
  ways. By including the permutation  $\alpha$ , we get totally $d$ sets of
  $(d-1)!$ permutations such that the permutations are not adjacent with each other
  with in each $d$ sets, but all $(d-1)!$ permutations  in each one of
  these $d$ sets are adjacent to the remaining sets of vectors. Thus
  they form a complete $d$-partite graph $K_d((d-1)!)$.
  Now the permutation $\alpha$  can vary in $n!/d!$ ways. Thus correspondingly
  we get $\frac{n!}{d!}$ copies of $K_d((d-1)!)$. i.e.,

\begin{center}
$G(S_n,\{d\})\cong \frac{n!}{d!} K_d((d-1)!)$.
\end{center}

 So for each elements of distance set $D$ we have
  $G (S_n,\{d_i\} )\cong \frac{n!}{d!} K_{d_i}((d_i-1)!)$

To get $G(S_n,D)$ we follow the same procedure as used in theorem
3.1.\hfill{$\Box$}\\

\noindent\textit{\textbf{Note 4.1: }}Let $ r_i ,s_i$
$(i=1,\ldots,k)$ be positive integers with $r_1 < \ldots <r_k$ and
$s_1 < \ldots < s_k$ then by using induction it is easy to show that
$r_1r_1!+ \ldots + r_k r_k!= s_1 s_1!+ \ldots + s_ks_k!$ if and only
if $r_i = s_i$, $\forall i=1,\ldots,k$.\\


\noindent\textit{ \textbf{Corollary 4.1:} For $(S_n,\rho)$, where
$\rho$ denotes the RT metric, suppose $D=\{d_1,\ldots
,d_k\}\subseteq $ dist $(S_n)$ with $ 2 \leq d_1 < \ldots < d_k \leq
n $ then the following holds.
\begin{itemize}
\item[(1)] $G(S_n,D)$ is regular of degree
 $\displaystyle \sum _{i=1}^{k}(d_i-1)(d_i-1)!$.
\item[(2)] $G(S_n,D)$ is connected if and only if $n \in D$.
\item[(3)] The components of $G(S_n,D)$ and $G(S_m, D_1)$ are
isomorphic if and only if $D = D_1$.
\item[(4)] $G(S_n,D)$ is a subgraph of $G(Z_n^n,D)$.
\item[(5)] $\chi \left( G (S_n,D )\right)= d_1d_2\ldots d_k$.
\end{itemize}}


\noindent\textit{\textbf{Proof :}} Proof of (1),(2),(4),(5) follows
directly form theorem 4.2. Proof of (3) follows form note 4.1 and
theorem 4.2.

 \end{document}